\documentclass[a4paper,12pt,final]{amsart}
\usepackage{times,a4wide,mathrsfs,amssymb,amsmath,amsthm,enumerate,xypic,tikzsymbols,dsfont}

\newcommand{\C}{\mathbb{C}}

\newcommand{\QQ}{\mathbb{Q}}
\newcommand{\NN}{\mathbb{N}}
\newcommand{\PP}{\mathbb{P}}

\newcommand{\OO}{\mathcal O}

\newcommand{\XX}{\mathcal X}
\newcommand{\YY}{\mathcal Y}

\newcommand{\MM}{\mathcal M}

\newcommand{\coker}{\hbox{Coker}}

\newcommand{\rom}{\romannumeral}

\newcommand{\one}{\mathds{1}}

\DeclareMathOperator{\aut}{Aut}

\DeclareMathOperator{\ide}{id}

\DeclareMathOperator{\ima}{Im}

\DeclareMathOperator{\sym}{Sym}
\DeclareMathOperator{\Gr}{Gr}
\DeclareMathOperator{\LGr}{LGr}
\DeclareMathOperator{\OGr}{OGr}

\newtheorem{theorem}{Theorem}[section]

\newtheorem{lemma}[theorem]{Lemma}

\newtheorem{proposition}[theorem]{Proposition}

\newtheorem{remark}[theorem]{Remark}
\newtheorem{definition}[theorem]{Definition}
\newtheorem{convention}{Conventions}
\newtheorem{question}[theorem]{Question}

\newtheorem{nonumbering}{Theorem}

\newtheorem{nonumberingc}{Corollary}

\newtheorem{nonumberingp}{Proposition}

\newtheorem{nonumberingt}{Acknowledgements}

\begin{document}

\author[Robert Laterveer]
{Robert Laterveer}

\address{Institut de Recherche Math\'ematique Avanc\'ee,
CNRS -- Universit\'e 
de Strasbourg,\
7 Rue Ren\'e Des\-car\-tes, 67084 Strasbourg CEDEX,
FRANCE.}
\email{robert.laterveer@math.unistra.fr}

\title[Some more Fano threefolds with an MCK decomposition]{Some more Fano threefolds with a multiplicative Chow--K\"unneth decomposition}

\begin{abstract} We exhibit several families of Fano threefolds with a multiplicative Chow--K\"unneth decomposition, in the sense of Shen--Vial. As a 
consequence, a certain tautological subring of the Chow ring of powers of these threefolds injects into cohomology. As a by-product of the argument, we observe that double covers of projective spaces admit a multiplicative Chow--K\"unneth decomposition.
 \end{abstract}

\thanks{Supported by ANR grant ANR-20-CE40-0023.}


\maketitle

\section{Introduction}

Given a smooth projective variety $Y$ over $\C$, let $A^i(Y):=CH^i(Y)_{\QQ}$ denote the Chow groups of $Y$ (i.e. the groups of codimension $i$ algebraic cycles on $Y$ with $\QQ$-coefficients, modulo rational equivalence). The intersection product defines a ring structure on $A^\ast(Y)=\bigoplus_i A^i(Y)$, the {\em Chow ring\/} of $Y$ \cite{F}. 

In the special case of K3 surfaces, this ring structure has remarkable properties:

\begin{theorem}[Beauville--Voisin \cite{BV}]\label{bv} Let $S$ be a projective K3 surface. 
The $\QQ$-subalgebra
  \[   \bigl\langle  A^1(S), c_j(S) \bigr\rangle\ \ \ \subset\ A^\ast(S) \]
  injects into cohomology under the cycle class map.
  \end{theorem}

\begin{theorem}[Voisin \cite{V17}, Yin \cite{Yin}]\label{vy} Let $S$ be a projective K3 surface, and $m\in\NN$. The $\QQ$-subalgebra
  \[   R^\ast(S^m):=\bigl\langle A^1(S), \Delta_S\bigr\rangle\ \ \ \subset\ A^\ast(S^m) \]
  (generated by pullbacks of divisors and pullbacks of the diagonal $\Delta_S\subset S\times S$)
  injects into cohomology under the cycle class map for all $m\le 2\dim H^2_{tr}(S,\QQ)+1$ (where $H^2_{tr}(S,\QQ)$ denotes the transcendental part of cohomology). Moreover, $R^\ast(S^m)$ injects into cohomology for all $m\in\NN$ if and only if $S$ is Kimura finite-dimensional.
  \end{theorem}

The Chow ring of abelian varieties also has an interesting property: there is a multiplicative splitting, defined by the Fourier transform \cite{Beau}.

Motivated by the particular behaviour of K3 surfaces and abelian varieties, Beauville \cite{Beau3} has conjectured that for certain special varieties, the Chow ring should admit a multiplicative splitting. In the wake of Beauville's ``splitting property conjecture'', Shen--Vial \cite{SV} have introduced the concept of {\em multiplicative Chow--K\"unneth decomposition\/}
(we will abbreviate this to ``MCK decomposition''). With the concept of MCK decomposition, it is possible to make concrete sense of this elusive ``splitting property conjecture'' of Beauville.

It is hard to understand precisely which varieties admit an MCK decomposition. To give an idea of what is known: hyperelliptic curves have an MCK decomposition \cite[Example 8.16]{SV}, but the very general curve of genus $\ge 3$ does not have an MCK decomposition \cite[Example 2.3]{FLV2}; K3 surfaces have an MCK decomposition, but certain high degree surfaces in $\PP^3$ do not have an MCK decomposition (cf. the examples given in \cite{OG}, cf. also section \ref{ss:mck} below). 

In this note, we will focus on Fano threefolds and ask the following question:

\begin{question}\label{ques} Let $X$ be a Fano threefold with Picard number 1. Does $X$ admit an MCK decomposition ?
\end{question} 

The restriction on the Picard number is necessary to rule out a counterexample of Beauville \cite[Examples 9.1.5]{Beau3}. The answer to Question \ref{ques} is affirmative for cubic threefolds \cite{Diaz}, \cite{FLV2}, for intersections of 2 quadrics \cite{2q}, for intersections of a quadric and a cubic \cite{55}, and for prime Fano threefolds of genus 8 \cite{g8} and of genus 10 \cite{g10}.

The main result of this paper answers Question \ref{ques} for several more families of Fano threefolds:

\begin{nonumbering}[=Theorem \ref{main}] The following smooth Fano threefolds have a multiplicative Chow--K\"unneth decomposition:

\begin{itemize}

\item hypersurfaces of weighted degree 6 in weighted projective space $\PP(1^3,2,3)$;

\item quartic double solids;

\item sextic double solids;

\item double covers of a quadric in $\PP^4$ branched along the intersection with a quartic;

\item special Gushel--Mukai threefolds.
\end{itemize}
\end{nonumbering}

In Table 1 (at the end of this paper), we have listed all Fano threefolds of Picard number 1 and what is known about MCK for them.

To prove Theorem \ref{main}, we provide a general criterion (Proposition \ref{crit}), that may be useful in other situations. For example, using this criterion we also prove the following:

\begin{nonumberingp}[=Proposition \ref{double}] Let $X$ be a smooth projective variety such that $X\to\PP^n$ is a double cover ramified along a smooth divisor $D\subset\PP^n$ of degree $d>n$. Then $X$ admits an MCK decomposition.
\end{nonumberingp}

As a consequence of Theorem \ref{main}, we obtain an injectivity result similar to Theorem \ref{vy}:

\begin{nonumberingc}[cf. Theorem \ref{taut}] Let $Y$ be a Fano threefold as in Theorem \ref{main}, and $m\in\NN$. Let
  \[ R^\ast(Y^m):=\bigl\langle h, \Delta_Y\bigr\rangle\ \subset\ \ \ A^\ast(Y^m)   \]
  be the $\QQ$-subalgebra generated by pullbacks of the polarization $h\in A^1(Y)$ and pullbacks of the diagonal $\Delta_Y\in A^3(Y\times Y)$. 
  The cycle class map induces injections
   \[ R^\ast(Y^m)\ \hookrightarrow\ H^\ast(Y^m,\QQ)\ \ \ \hbox{for\ all\ }m\in\NN\ .\]
   \end{nonumberingc}


 \vskip0.6cm

\begin{convention} In this paper, the word {\sl variety\/} will refer to a reduced irreducible scheme of finite type over $\C$. A {\sl subvariety\/} is a (possibly reducible) reduced subscheme which is equidimensional. 

{\bf All Chow groups will be with rational coefficients}: we will denote by $A_j(Y)$ the Chow group of $j$-dimensional cycles on $Y$ with $\QQ$-coefficients; for $Y$ smooth of dimension $n$ the notations $A_j(Y)$ and $A^{n-j}(Y)$ are used interchangeably. 
The notation $A^j_{hom}(Y)$ will be used to indicate the subgroup of homologically trivial cycles.
For a morphism $f\colon X\to Y$, we will write $\Gamma_f\in A_\ast(X\times Y)$ for the graph of $f$.

The contravariant category of Chow motives (i.e., pure motives with respect to rational equivalence as in \cite{Sc}, \cite{MNP}) will be denoted 
$\MM_{\rm rat}$.
\end{convention}

\section{MCK decomposition}
\label{ss:mck}

\begin{definition}[Murre \cite{Mur}] Let $X$ be a smooth projective variety of dimension $n$. We say that $X$ has a {\em CK decomposition\/} if there exists a decomposition of the diagonal
   \[ \Delta_X= \pi^0_X+ \pi^1_X+\cdots +\pi_X^{2n}\ \ \ \hbox{in}\ A^n(X\times X)\ ,\]
  such that the $\pi^i_X$ are mutually orthogonal idempotents and $(\pi_X^i)_\ast H^\ast(X,\QQ)= H^i(X,\QQ)$.
  
  (NB: ``CK decomposition'' is shorthand for ``Chow--K\"unneth decomposition''.)
\end{definition}

\begin{remark} Murre has conjectured that any smooth projective variety should have a CK decomposition \cite{Mur}, \cite{J4}. 
\end{remark}

\begin{definition}[Shen--Vial \cite{SV}] Let $X$ be a smooth projective variety of dimension $n$, and let $\Delta_X^{sm}\in A^{2n}(X\times X\times X)$ denote the class of the small diagonal
  \[ \Delta_X^{sm}:=\bigl\{ (x,x,x)\ \vert\ x\in X\bigr\}\ \subset\ X\times X\times X\ .\]
  An {\em MCK decomposition\/} is defined as a CK decomposition $\{\pi_X^i\}$ of $X$ that is {\em multiplicative\/}, i.e. it satisfies
  \[ \pi_X^k\circ \Delta_X^{sm}\circ (\pi_X^i\times \pi_X^j)=0\ \ \ \hbox{in}\ A^{2n}(X\times X\times X)\ \ \ \hbox{for\ all\ }i+j\not=k\ .\]
  
 (NB: ``MCK decomposition'' is shorthand for ``multiplicative Chow--K\"unneth decomposition''.) 
  \end{definition}
  
  \begin{remark} The small diagonal (when considered as a correspondence from $X\times X$ to $X$) induces the {\em multiplication morphism\/}
    \[ \Delta_X^{sm}\colon\ \  h(X)\otimes h(X)\ \to\ h(X)\ \ \ \hbox{in}\ \MM_{\rm rat}\ .\]
 Let us assume $X$ has a CK decomposition
  \[ h(X)=\bigoplus_{i=0}^{2n} h^i(X)\ \ \ \hbox{in}\ \MM_{\rm rat}\ .\]
  By definition, this decomposition is multiplicative if for any $i,j$ the composition
  \[ h^i(X)\otimes h^j(X)\ \to\ h(X)\otimes h(X)\ \xrightarrow{\Delta_X^{sm}}\ h(X)\ \ \ \hbox{in}\ \MM_{\rm rat}\]
  factors through $h^{i+j}(X)$.
  
  If $X$ has an MCK decomposition, then setting
    \[ A^i_{(j)}(X):= (\pi_X^{2i-j})_\ast A^i(X) \ ,\]
    one obtains a bigraded ring structure on the Chow ring: that is, the intersection product sends $A^i_{(j)}(X)\otimes A^{i^\prime}_{(j^\prime)}(X) $ to  $A^{i+i^\prime}_{(j+j^\prime)}(X)$.
    
      It is conjectured that for any $X$ with an MCK decomposition, one has
    \[ A^i_{(j)}(X)\stackrel{??}{=}0\ \ \ \hbox{for}\ j<0\ ,\ \ \ A^i_{(0)}(X)\cap A^i_{hom}(X)\stackrel{??}{=}0\ ;\]
    this is related to Murre's conjectures B and D, that have been formulated for any CK decomposition \cite{Mur}.


For more background on the concept of MCK, and for examples of varieties with an MCK decomposition, we refer to \cite[Section 8]{SV}, as well as \cite{V6}, \cite{SV2}, \cite{FTV}, \cite{37},
  \cite{38}, \cite{39}, \cite{40}, \cite{44}, \cite{FLV2}, \cite{60}, \cite{55}, \cite{59}, \cite{NOY}.
    \end{remark}

 \section{A general criterion}
 
 We develop a general criterion for having an MCK. The criterion hinges on the {\em Franchetta property\/} for families of varieties, which is defined as follows:

  \begin{definition}\label{def} Let $\XX\to B$ be a smooth projective morphism, where $\XX, B$ are smooth quasi-projective varieties, and let us write $X_b$ for the fiber over $b\in B$. We say that $\XX\to B$ has the {\em Franchetta property in codimension $j$\/} if the following holds: for every $\Gamma\in A^j(\XX)$ such that the restriction $\Gamma\vert_{X_b}$ is homologically trivial for the very general $b\in B$, the restriction $\Gamma\vert_b$ is zero in $A^j(X_b)$ for all $b\in B$.
 
 We say that $\XX\to B$ has the {\em Franchetta property\/} if $\XX\to B$ has the Franchetta property in codimension $j$ for all $j$.
 \end{definition}
 
 This property is studied in \cite{PSY}, \cite{BL}, \cite{FLV}, \cite{FLV3}.
 
 \begin{definition} Given a family $\XX\to B$ as in Definition \ref{def}, we use the shorthand
   \[ GDA^j_B(X_b):=\ima\bigl( A^j(\XX)\to A^j(X_b)\bigr)\ \ \ \subset\ A^j(X_b) \ \]
   ($GDA^\ast()$ stands for the ``generically defined cycles'').
  \end{definition}
 
 The Franchetta property for $\XX\to B$ means that the generically defined cycles inject into cohomology.

%
%
%
%
%
%
%
 
 \begin{proposition}\label{crit} Let $\XX\to B$ be a family of smooth projective varieties of relative dimension $n$, with fiber $X_b$. Assume the following:
 
 (\rom1) the family $\XX\times_B \XX\to B$ has the Franchetta property;
  
 (\rom2) there exists a projective quotient variety $P$ (i.e. $P=P^\prime/G$ where $P^\prime$ is smooth projective and $G\subset\aut(P^\prime)$ is a finite cyclic group) with trivial Chow groups (i.e. $A^\ast_{hom}(P)=0$), such that $X_b\to P$ is a double cover with branch locus a smooth ample divisor, for all $b\in B$.
 
 Then $X_b$ admits an MCK decomposition, for all $b\in B$.
 \end{proposition}
 
 \begin{proof} 
 We have the following Lefschetz-type result in cohomology:
 
 \begin{lemma}\label{cor} Let $X_b\to P$ be as in the proposition.
 Then pullback
   \[ H^i(P,\QQ)\ \to\ H^i(X_b,\QQ) \]
   is an isomorphism for $i<n$, and injective for $i=n$.
   \end{lemma}
   
   \begin{proof} In case $P$ is smooth, this is a result of Cornalba \cite{Cor}. The general case is readily deduced from this: assume $P=P^\prime/G$ where $P^\prime$ is smooth projective and $G\subset\aut(P^\prime)$ is a finite cyclic group, and consider the fiber square
     \[ \begin{array}[c]{ccc}
             X_b^\prime & \to& X_b\\
             &&\\
             \downarrow&&\downarrow\\
             &&\\
             P^\prime&\to&\ P\ .\\
             \end{array}\]
     Cornalba's result applies to the double cover of the left-hand vertical arrow, and so pullback        
      \[ H^i(P^\prime,\QQ)\ \to\ H^i(X^\prime_b,\QQ) \]
   is an isomorphism for $i<n$, and injective for $i=n$.
   The $G$-action on $P^\prime$ lifts to $X^\prime_b$, and taking $G$-invariants we find that
    \[ H^i(P,\QQ)=H^i(P^\prime,\QQ)^G\ \to\ H^i(X^\prime_b,\QQ)^G=H^i(X_b,\QQ) \]
   is an isomorphism for $i<n$, and injective for $i=n$.
          \end{proof}

   Since $H^\ast(P,\QQ)$ is algebraic (this is a general fact for any variety with trivial Chow groups, cf. \cite{Kim2}), this implies that also
   $H^i(X_b,\QQ)$ is algebraic, for all $i\not=n$. More precisely, for $i\not= n$ odd, one has $H^i(X_b,\QQ)=0$ while for $i<n$ even, one has isomorphisms
     \[ A^{i/2}(P)\ \cong\ H^i(X_b,\QQ) \ ,\]
     induced by pullback. This implies that for $i<n$ the K\"unneth components $\pi^i_{X_b}$ are algebraic, and generically defined. 
     To define the K\"unneth components $\pi^i_{X_b}$ explicitly, let $p\colon X_b\to P$ denote the projection morphism, and let $\pi^i_P$ denote the (unique) CK decomposition of $P$.
     One can then define
       \[ \begin{split}   \pi_{X_b}^i&:= {1/ 2}\, {}^t  \Gamma_p\circ \pi^i_P\circ \Gamma_p\ \ \ \hbox{if}\ i<n\ ,\\
                                 \pi_{X_b}^i&:=\pi^{2n-i}_{X_b}\ \ \ \ \ \ \hbox{if}\ i>n\ ,\\
                                 \pi_{X_b}^{n, fix}&:= {1/ 2}\, {}^t  \Gamma_p\circ \pi^{n}_P\circ \Gamma_p\    ,\\
                                 \pi_{X_b}^{n, var}&:=\Delta_{X_b}-\sum_{j\not=n}\pi_{X_b}^j -\pi_{X_b}^{n,fix}\ ,\\
                                 \pi_{X_b}^{n}&:=    \pi_{X_b}^{n, fix}+    \pi_{X_b}^{n, var}\ \ \ \in\ A^n(X_b\times X_b)\ .\\
                                 \end{split}           \]
              (Note that $\pi^n_{X_b}=0$ in case $n$ is odd.)     
              The notation is meant to remind the reader that $\pi_{X_b}^{n, fix}$ and $\pi_{X_b}^{n, var}$ are projectors on the fixed part resp. the variable part of cohomology in degree $n$.

       These projectors define a generically defined CK decomposition for each $X_b$, i.e. all projectors are in $GDA^n_B(X_b\times X_b)$.
       This CK decomposition has the property that
       \begin{equation}\label{one} \begin{split} h^j(X_b)&:=(X_b,\pi^j_{X_b},0)=\oplus \one(\ast)\ \ \ \forall j\not=n\ ,\\
           h^{n,fix}(X_b)&:=(X_b,\pi_{X_b}^{n,fix},0)=\oplus \one(\ast)\ \ \hbox{in}\ \MM_{\rm rat}\ .\\
           \end{split} \end{equation}
       
   Let us now proceed to verify that this CK decomposition is MCK. What we need to check is the vanishing
          \[ \pi_{X_b}^k\circ \Delta_{X_b}^{sm}\circ (\pi_{X_b}^i\times \pi_{X_b}^j)=0\ \ \ \hbox{in}\ A^{2n}(X_b\times X_b\times X_b)\ \ \ \hbox{for\ all\ }i+j\not=k\ .\]
     
  First, let us assume that at least one of the 3 integers $(i,j,k)$ is different from $n$, and $i+j\not=k$. In this case, we have that
    \[ \begin{split} \pi_{X_b}^k\circ \Delta_{X_b}^{sm}\circ (\pi_{X_b}^i\times \pi_{X_b}^j)&= ({}^t     \pi_{X_b}^i\times {}^t \pi_{X_b}^j\times \pi^k_{X_b})_\ast \Delta^{sm}_{X_b}\\
                                                           & = ({}   \pi_{X_b}^{2n-i}\times {}\pi_{X_b}^{2n-j}\times \pi^k_{X_b})_\ast \Delta^{sm}_{X_b}\\  
                                                           & \hookrightarrow \bigoplus A^\ast(X_b\times X_b)\ .\\
                                                           \end{split}\]
                                                           Here the first equality is an application of Lieberman's lemma \cite[Lemma 2.1.3]{MNP}, and the inclusion follows from property \eqref{one}. The resulting cycle in $\bigoplus A^\ast(X_b\times X_b)$ is generically defined (since the $\pi^\ast_{X_b}$ and $\Delta_{X_b}^{sm}$ are) and homologically trivial (since $i+j\not=k$). By assumption (\rom1), the resulting cycle in $\bigoplus A^\ast(X_b\times X_b)$ is rationally trivial, and so        
                                                          \[ \pi_{X_b}^k\circ \Delta_{X_b}^{sm}\circ (\pi_{X_b}^i\times \pi_{X_b}^j)=0\ \ \ \hbox{in}\ A^{2n}(X_b\times X_b\times X_b)\ ,\]
                                                 as desired.
                                                 
                                      It remains to treat the case $i=j=k=n$. The decomposition $ \pi_{X_b}^{n}:=    \pi_{X_b}^{n, fix}+    \pi_{X_b}^{n, var}$ induces a decomposition
                                      \[  \begin{split} \pi_{X_b}^n\circ \Delta_{X_b}^{sm}\circ (\pi_{X_b}^n\times \pi_{X_b}^n)=& \pi_{X_b}^{n,fix}\circ \Delta_{X_b}^{sm}\circ (\pi_{X_b}^{n,fix}\times 
                                      \pi_{X_b}^{n,fix}) \\
                                       &+    \pi_{X_b}^{n,fix}\circ \Delta_{X_b}^{sm}\circ (\pi_{X_b}^{n,fix}\times \pi_{X_b}^{n,var})  \\
                                       &+ \cdots\ \cdots\\
                                       &+ \pi_{X_b}^{n,var}\circ \Delta_{X_b}^{sm}\circ (\pi_{X_b}^{n,var}\times \pi_{X_b}^{n,var})
                                      \ \ \ \hbox{in}\ A^{2n}(X_b\times X_b\times X_b)\ .\\
                                      \end{split}
                                      \]   
     Using property \eqref{one} and the Franchetta property for $X_b\times X_b$, all summands containing $\pi_{X_b}^{n,fix}$ vanish. One is left with the last term. To deal with the last term,
     we observe that the covering involution $\iota\in\aut(X_b)$ of the double cover $p\colon X_b  \to P$ induces a splitting of the motive
       \[ \begin{split} h(X_b)=& h(X_b)^+\oplus h(X_b)^-\\
                                         :=&(X_b, 1/2\, (\Delta_{X_b}+\Gamma_\iota),0)\oplus (X_b, 1/2\, (\Delta_{X_b}-\Gamma_\iota),0)       \ \ \ \hbox{in}\ \MM_{\rm rat}\ ,
                                         \end{split}\]
                                         where $\Gamma_\iota$ denotes the graph of the involution $\iota$.
                      Moreover, there is equality
                           \[  h^{n,var}(X_b)=h(X_b)^-\ \ \ \hbox{in}\ \MM_{\rm rat}\ .\]         
            But the intersection product map
            \[  h(X_b)^-\otimes h(X_b)^-\ \ \xrightarrow{\Delta^{sm}_{X_b}}\ h(X_b) \]
            factors over $h(X_b)^+$, as is readily seen (cf. Lemma \ref{+-} below), which is saying exactly that
            \[     \pi_{X_b}^{n,var}\circ \Delta_{X_b}^{sm}\circ (\pi_{X_b}^{n,var}\times \pi_{X_b}^{n,var})=0
                                      \ \ \ \hbox{in}\ A^{2n}(X_b\times X_b\times X_b)\ .\]
           
           This closes the proof, modulo the following lemma (which is probably well-known, but we include a proof for completeness):
           
           \begin{lemma}\label{+-} Let $ X\to P$ be a double cover, where $X$ and $P$ are quotient varieties, and let $\iota\in\aut(X)$ be the covering involution. Let
             \[ h(X)^+:= (X, 1/2\, (\Delta_{X}+\Gamma_\iota,0)\  , \ \ \ h(X)^-:=(X, 1/2\, (\Delta_{X}-\Gamma_\iota),0)\ \ \ \hbox{in}\ \MM_{\rm rat}\ .\]
         The map of motives
          \[  h(X)^-\otimes h(X)^-\ \ \xrightarrow{\Delta^{sm}_{X}}\ h(X) \]
            factors over $h(X)^+$.             
                  \end{lemma}    
                  
                To prove the lemma, let $\iota\in\aut(X)$ denote the covering involution. 
                The motive $h(X)^-$ is defined by the projector
                \[ \Delta_X^-:= 1/2\, (\Delta_X-\Gamma_\iota)\ \ \ \in\ A^n(X\times X)\ .\]
                Plugging this in and developing, it follows that
                \[ \begin{split}
                        \Delta_X^-\circ \Delta^{sm}_X\circ (\Delta_X^-\times\Delta_X^-) &=1/8\, (\Delta_X-\Gamma_\iota)\circ \Delta_X^{sm}\circ (\Delta_{X\times X} -\Delta_X\times\Gamma_\iota - \Gamma_\iota\times\Delta_X + \Gamma_\iota\times\Gamma_\iota)\\
                        &= 1/8\, \Bigl( \Delta_X\circ\Delta_X^{sm}\circ(\Delta_X\times\Delta_X) + \cdots - \Gamma_\iota\circ\Delta_X^{sm}\circ(\Gamma_\iota\times\Gamma_\iota)\Bigr)\\
                        &=1/8\, \Bigl( \Delta_X^{sm}\\ 
                                 &   \ \ \ \  \ \ \ \  - (\ide\times\ide\times\iota)_\ast(\Delta_X^{sm}) -  (\ide\times\iota\times\ide)_\ast(\Delta_X^{sm}) - (\iota\times\ide\times\ide)_\ast(\Delta_X^{sm}) \\
                        & \ \ \ \ \ \ \ \ +  (\ide\times\iota\times\iota)_\ast(\Delta_X^{sm}) +  (\iota\times\ide\times\iota)_\ast(\Delta_X^{sm})   + (\iota\times\iota\times\ide)_\ast(\Delta_X^{sm})\\
                        & \ \ \ \ \ \  \ \ -
                           (\iota\times\iota\times\iota)_\ast(\Delta_X^{sm}) \Bigr)\ \ \ \hbox{in}\ A^{2n}(X\times X\times X)\ .\\
                         \end{split} \]
                         Here the last equality is by virtue of Lieberman's lemma \cite[Lemma 2.1.3]{MNP}.
                         However, we have equality
                         \[ \Delta_X^{sm}=\{(x,x,x)\ \vert\ x\in X\} = (\iota\times\iota\times\iota)_\ast (\Delta_X^{sm})\ \ \ \hbox{in}\ A^{2n}(X\times X\times X)\ ,\]
                         and so the sum of the first and last summand vanish. Likewise, we have equality
                         \[ (\ide\times\iota\times\iota)_\ast(\Delta_X^{sm}) = (\ide\times\iota\times\iota)_\ast (\iota\times\iota\times\iota)_\ast(\Delta_X^{sm}) = (\iota\times\ide\times\ide)_\ast(\Delta_X^{sm})\ \ \ \hbox{in}\ A^{2n}(X\times X\times X)\ ,\]
                         and so the other summands cancel each other pairwise. This proves the lemma.
                         \end{proof}

%

As a first application of our general criterion, we now proceed to show the following:

\begin{proposition}\label{double} Let $X$ be a smooth projective variety such that $X\to\PP^n$ is a double cover ramified along a smooth divisor $D\subset\PP^n$, and assume either
$\dim H^n(X,\QQ)>1$, or $D$ has degree $d>n$. Then $X$ admits an MCK decomposition.
\end{proposition}

\begin{proof} Double covers $X$ as in the proposition are exactly the smooth hypersurfaces of degree $2d$ in the weighted projective space $\PP:=\PP(1^{n+1},d)$, where $2d:=\deg D$.
Let 
  \[ B\subset \bar{B}:=\PP H^0(\PP,\OO_\PP(2d)) \]
  denote the Zariski open parametrizing smooth hypersurfaces, and let 
  \[ B\times\PP    \  \supset \  \XX\ \to\ B \]
  denote the universal family. In view of Proposition \ref{crit}, it suffices to check that the family $\XX\times_B \XX\to B$ has the Franchetta property.
  
  To this end, we remark that the line bundle $\OO_\PP(2d)$ is very ample (cf. Lemma \ref{va} below),
  which means that the set-up verifies condition $(\ast_2)$ of \cite[Definition 2.5]{FLV2}. 
  An application of the stratified projective bundle argument 
  \cite[Proposition 2.6]{FLV2} then implies that
    \begin{equation}\label{gd0}   GDA^\ast_B(X_b\times X_b) = \bigl\langle (p_i)^\ast(h) ,\Delta_{X_b}\bigr\rangle\ ,\end{equation}
    where we write $h\in A^1(X_b)$ for the hyperplane class. The excess intersection formula \cite[Theorem 6.3]{F} gives an equality
    \[ \Delta_{X_b}\cdot (p_i)^\ast(h) = 2d\, \sum_j (p_1)^\ast(h^j)\cdot (p_2)^\ast(h^{n+1-j})\ \ \ \hbox{in}\ A^{n+1}(X_b\times X_b)\ ,\]
    and so equality \eqref{gd0} reduces to the equality
    \[  GDA^\ast_B(X_b\times X_b) = \bigl\langle (p_1)^\ast(h) , (p_2)^\ast(h)\bigr\rangle \oplus \QQ[\Delta_{X_b}]\ .\]
    The ``decomposable part''  $\langle (p_1)^\ast(h) , (p_2)^\ast(h)\rangle$  injects into cohomology, because of the K\"unneth formula for $H^\ast(X_b\times X_b,\QQ)$. 
    The class of the diagonal in cohomology is linearly independent from the decomposable part: indeed, if the diagonal were decomposable it would act as zero on the primitive cohomology
    \[ H^n_{\rm prim}(X_b,\QQ):=\coker\bigl( H^n(\PP^n,\QQ)\to H^n(X_b,\QQ)\bigr) \ .\]
    But the assumption $\dim H^n(X_b,\QQ)>1$ is equivalent to having $H^n_{\rm prim}(X_b,\QQ)\not=0$. This proves the Franchetta property for $\XX\times_B \XX\to B$, and closes the proof.              
    
    The case $d>n$ is a special case where $H^n_{\rm prim}(X_b,\QQ)\not=0$, because it is known that    
     the geometric genus of $X_b$ is
     \[ p_g(X_b)= {d-1\choose n} \]
      \cite[Section 3.5.4]{Dol}.

      It remains to prove the following, which we have used above:
      
      \begin{lemma}\label{va} Let $\PP:=\PP(1^{n+1},d)$.  The sheaf $\OO_\PP(d)$ is locally free and very ample.
      \end{lemma}
      
      The assertion about the sheaf being locally free is just because $d$ is a multiple of the weights of $\PP$ (cf. \cite[Remarques 1.8]{Del}). As for the very ampleness, we apply Delorme's criterion \cite[Proposition 2.3(\rom3)]{Del} (cf. also \cite[Theorem 4.B.7]{BR}). To prove very ampleness of $\OO_\PP(d)$, we need to prove that the integer $E$ as defined in \cite{Del} and \cite{BR} is equal to 0.

Let us write $x_0, \ldots, x_n,y$ for the weighted homogeneous coefficients of $\PP$, where $x_j$ and $y$ have weight $1$ resp. $d$. It is readily seen that every monomial in $x_j,y$ of (weighted) degree $m+dk$ (where $m$ is a positive multiple of $d$, and $k$ is any positive integer) is divisible by a monomial of (weighted) degree $dk$. This means that the integer $E$ defined in loc. cit. is 0, and so \cite[Proposition 2.3(\rom3)]{Del} implies the very ampleness of $\OO_\PP(d)$. 

This proves the lemma, and ends the proof of the proposition.
                   \end{proof}

  Here is another sample application of our general criterion:    

\begin{proposition}\label{123} Let $X\subset\PP(1^n, 2,3)$ be a smooth hypersurface of (weighted) degree 6. Assume $\dim H^n(X,\QQ)>1$. Then $X$ has an MCK decomposition.
\end{proposition}

\begin{proof} The varieties $X$ as in the proposition are exactly the smooth double covers of $\PP:=\PP(1^n,2)$ branched along a (weighted) degree 6 divisor (cf. \cite[Remark 2.3]{KP} and for $n=3$ also \cite[Theorem 4.2]{Isk}).
Let $\XX\to B$ denote the family of such double covers. We are going to check that the family $\XX\times_B \XX\to B$ has the Franchetta property. Proposition \ref{123} is then a special case of our general criterion Proposition \ref{crit}.

Let $\bar{\XX}\to\bar{B}\cong\PP^r$ denote the universal family of all (possibly singular) hypersurfaces of weighted degree 6 in $\PP$. The line bundle $\OO_\PP(6)$ is very ample
(cf. Lemma \ref{very} below), and so the projection
  \[  \bar{\XX}\times_{\bar{B}}\bar{\XX}\ \to\ \PP\times\PP \]
  has the structure of a stratified projective bundle (with strata the diagonal $\Delta_\PP$ and its complement). One can thus use the stratified projective bundle argument \cite[Proposition 2.6]{FLV2} to deduce the identity
    \[ \begin{split} GDA^\ast_B(X\times X)&= \bigl\langle  (p_i)^\ast GDA^\ast_B(X), \Delta_X\bigr\rangle\\
                                                                        &= \bigl\langle  (p_i)^\ast (h), \Delta_X\bigr\rangle\\
                                                                        \end{split}\]
                                                                        (here, $h\in A^1(X)$ denotes the restriction to $X$ of an ample generator of $A^1(\PP)\cong\QQ$).
                                                                        
Since $X\subset\PP$ is a hypersurface, the excess intersection formula gives
  \[   \Delta_X\cdot (p_i)^\ast(h) = \Delta_\PP\vert_X \ \ \in\ \bigl\langle  (p_i)^\ast (h)\bigr\rangle\ .\]
The
above identification thus simplifies to
  \[   GDA^\ast_B(X\times X)  = \bigl\langle  (p_i)^\ast (h)\bigr\rangle\ \oplus\ \QQ[\Delta_X]\ .\]
  
The assumption that $\dim H^n(X,\QQ)>1$ implies that the diagonal $\Delta_X$ is linearly independent in cohomology from the decomposable classes $ \bigl\langle  (p_i)^\ast (h)\bigr\rangle$
(indeed, the decomposable classes act as zero on the primitive cohomology of $X$, while the diagonal acts as the identity). This shows that $ GDA^\ast_B(X\times X)$ injects into cohomology, as requested.

\begin{lemma}\label{very} Let $\PP:=\PP(1^n,2,3)$. The sheaf $\OO_\PP(6)$ is (locally free and) very ample.
\end{lemma}

The assertion about the sheaf being locally free is just because $6$ is a multiple of all the weights (cf. \cite[Remarques 1.8]{Del}). As for the very ampleness, we apply Delorme's criterion \cite[Proposition 2.3(\rom3)]{Del} (cf. also \cite[Theorem 4.B.7]{BR}). To prove very ampleness of $\OO_\PP(6)$, we need to prove that the integer $E$ defined in \cite{Del} and \cite{BR} is equal to 0.

Let us write $x_1, \ldots, y,z$ for the weighted homogeneous coefficients of $\PP$, where $y$ and $z$ have weight $2$ resp. $3$. We need to check that every monomial in $x_j,y,z$ of (weighted) degree $6+6k$ is divisible by a monomial of (weighted) degree $6k$ (if this is the case, then $E=0$ and \cite[Proposition 2.3(\rom3)]{Del} implies the very ampleness of $\OO_\PP(6)$). In case the monomial contains $z^2$, it is divisible by $z^2$ and so the condition is satisfied. Assume now the monomial contains only one $z$. In case the monomial contains $y^3$ it is divisible by $y^3$. Next, if the monomial contains $y$ (or $y^2$) it is divisible by $zyx_j$ (for some $j$) and so the condition is satisfied. A monomial in $z$ and $x_j$  obviously satisfies the condition. Finally, monomials in $x_j$ satisfy the condition.

This proves the lemma, and ends the proof of the proposition.
\end{proof}

 \section{Main result}
 
 \begin{theorem}\label{main} The following Fano threefolds admit an MCK decomposition:
 
 (\rom1) hypersurfaces of weighted degree 6 in weighted projective space $\PP(1^3,2,3)$;

(\rom2) quartic double solids;

(\rom3) sextic double solids;

(\rom4) double covers of a quadric in $\PP^4$ branched along the intersection with a quartic;

(\rom5) special Gushel--Mukai threefolds. 
  \end{theorem}
 
 \begin{proof} The cases (\rom2) and (\rom3) are immediate applications of Proposition \ref{double}. The case (\rom1) is a special case of Proposition \ref{123}.

 Before proving case (\rom4), let us first state a preparatory lemma:
 
   \begin{lemma}\label{lb} Let $Z\subset\PP:=\PP(1^5,2)$ be a smooth weighted hypersurface of degree 2. Then
     \[ \Delta_Z={1\over 2}\, \sum_{j=0}^4 h^j\times h^{4-j}\ \ \ \hbox{in}\ A^4(Z\times Z)\ .\]
   \end{lemma}
   
   \begin{proof} $Z$ is a quotient of a non-singular quadric in $\PP^5$ and so $Z$ has trivial Chow groups (i.e. $A^\ast_{hom}(Z)=0$). Using \cite[4.4.2]{Dol}, one can compute the Betti numbers of $Z$ and one finds that they are the same as those of projective space $\PP^4$. This means that there is a cohomological decomposition of the diagonal
     \[ \Delta_Z={1\over 2}\, \sum_{j=0}^4 h^j\times h^{4-j}\ \ \ \hbox{in}\ H^8(Z\times Z,\QQ)\ .\]   
     Since $Z$ (and hence also $Z\times Z$) has trivial Chow groups, the same decomposition holds modulo rational equivalence, proving the lemma.
      \end{proof}

Now, to prove case (iv) of Theorem \ref{main}, we apply our general criterion Proposition \ref{crit}.
 Let $\PP:=\PP(1^5,2)$, and let $\YY\to B$ be the universal family of smooth dimensionally transverse complete intersections of $\OO_\PP(2)\oplus\OO_\PP(4)$, where the base $B$ is a Zariski open 
   \[ B\subset\bar{B}:=\PP H^0(\PP,\OO_\PP(2)\oplus\OO_\PP(4))\ .\] 
   It follows from Lemma \ref{va} that $\OO_\PP(2)$ and $\OO_\PP(4)$ are very ample line bundles on $\PP$, and so $\bar{\YY}\times_{\bar{B}}\bar{\YY}\to\PP\times\PP$ is a stratified projective bundle with strata $\Delta_\PP$ and its complement.
 The usual stratified projective bundle argument \cite[Proposition 2.6]{FLV2} applies, and we find that
        \[ \begin{split} GDA^\ast_B(Y\times Y)&= \bigl\langle  (p_i)^\ast GDA^\ast_B(Y), \Delta_Y\bigr\rangle\\
                                                                        &= \bigl\langle  (p_i)^\ast (h), \Delta_Y\bigr\rangle\\
                                                                        \end{split}\]
                                                                        (here, $h\in A^1(Y)$ denotes the restriction to $Y$ of an ample generator of $A^1(\PP)\cong\QQ$).
Let $Y=Z\cap Z^\prime$, where $Z$ and $Z^\prime\subset\PP$ are hypersurfaces of (weighted) degree 2 and 4. Up to shrinking $B$, we may assume the hypersurface $Z$ is smooth.                                                                        
Since $Y\subset Z$ is a divisor, the excess intersection formula gives
  \[   \Delta_Y\cdot (p_i)^\ast(h) = \Delta_Z\vert_Y\ \ \ \hbox{in}\  A^4(Y\times Y)\ .\]
 Using Lemma \ref{lb}, it follows that
   \[  \Delta_Y\cdot (p_i)^\ast(h) \ \ \in\  \bigl\langle  (p_i)^\ast (h)\bigr\rangle\ .\]   
   The above identification thus simplifies to
  \[   GDA^\ast_B(Y\times Y)  = \bigl\langle  (p_i)^\ast (h)\bigr\rangle\ \oplus\ \QQ[\Delta_Y]\ .\]
   As before, the fact that the diagonal $\Delta_Y$ is linearly independent from the decomposable correspondences in cohomology now shows that
   \[   GDA^\ast_B(Y\times Y)\ \to\ H^\ast(Y\times Y,\QQ) \]
   is injective, and so $Y$ verifies the hypotheses of Proposition \ref{crit}.
      
  The argument for case (\rom5) is similar to that of (\rom4). 
  First, in view of the spread argument \cite[Lemma 3.2]{Vo}, it suffices to establish an MCK decomposition for the {\em generic\/} special Gushel--Mukai threefold $Y$. 
  Thus we may assume that there exists $P\subset \Gr(2,5)$, a smooth complete intersection of Pl\"ucker hyperplanes, and a double cover $p\colon Y\to P$ branched along a smooth Gushel--Mukai surface. We now consider the family $\YY\to B$ of all double covers of $P$ branched along smooth Gushel--Mukai surfaces (so $B\subset\bar{B}$ is a Zariski open in the projectivized space of quadratic sections of the cone over $P$), and we apply our general criterion Proposition \ref{crit} to this family. 
  
  \begin{lemma}\label{FY} Let $\YY\to B$ be the family of double covers of $P$ branched along smooth Gushel--Mukai surfaces. The family $\YY\to B$ has the Franchetta property.
  \end{lemma}
  
  \begin{proof} We consider the family $\bar{\YY}\to\bar{B}$ with the projection to the cone $C$ over $P$. This is a projective bundle, and so for any fiber $Y=Y_b$ with $b\in B$ we have
    \[  GDA^\ast_B(Y)=\ima\bigl( A^\ast(C)\to A^\ast(Y)\bigr)\ .\]
    The condition $b\in B$ means exactly that $Y$ avoids the summit of the cone $C$, and so (writing $C^\circ\subset C$ for the complement of the summit of the cone) we have
    \begin{equation}\label{gd}  GDA^\ast_B(Y)=\ima\bigl( A^\ast(C^\circ)\to A^\ast(Y)\bigr)\ .
                                                  \end{equation}
                                                  But $C^\circ\to P$ is an affine bundle, and 
                                                  \[A^\ast(P)=\ima\bigl(A^\ast(\Gr(2,5))\to A^\ast(P)\bigr) = \bigl\langle  h\bigr\rangle\ ,\] where $h$ denotes the restriction to $P$ of a Pl\"ucker hyperplane  (this follows from \cite[Theorem 3.17]{56}, or alternatively from the fact that the derived category of $P$ has a full exceptional collection of length 4 \cite{Orl}).
                                                  Thus, \eqref{gd} reduces to
                                                  \[  GDA^\ast_B(Y)=   \bigl\langle  h\bigr\rangle      \ .\]
                                                  This proves the Franchetta property for $Y$.
                                                    \end{proof}
                                                    
                                          \begin{lemma}\label{FY2} Let $\YY\to B$ be as in Lemma \ref{FY}. The family $\YY\times_B \YY\to B$ has the Franchetta property.
  \end{lemma}
  
  \begin{proof} Let us consider the family $\bar{\YY}\times_{\bar{B}} \bar{\YY}\to\bar{B}$ with the projection to $C\times C$. This is a stratified projective bundle, with strata $\Delta_C$ and its complement. Thus, the stratified projective bundle argument \cite[Proposition 2.6]{FLV2} implies that
  \[ GDA^\ast_B(Y\times Y)=\bigl\langle \ima\bigl(A^\ast(C^\circ\times C^\circ)\to A^\ast(Y\times Y)\bigr), \Delta_Y\bigr\rangle\ .\]
  Since $A^\ast(C^\circ)=\ima\bigl( A^\ast(\Gr(2,5))\to A^\ast(C^\circ)$, we find that
  \[   GDA^\ast_B(Y\times Y)=\bigl\langle \ima\bigl(A^\ast(\Gr(2,5)\times \Gr(2,5))\to A^\ast(Y\times Y)\bigr), \Delta_Y\bigr\rangle\ .\]  
  But $A^\ast(\Gr(2,5)\times \Gr(2,5))=A^\ast(\Gr(2,5))\otimes A^\ast(\Gr(2,5))  $ since the Grassmannian has trivial Chow groups, and so
  \[   \begin{split} GDA^\ast_B(Y\times Y)&=\bigl\langle   GD_B(Y), \Delta_Y \bigr\rangle \\  
  &= \bigl\langle  h, \Delta_Y \bigr\rangle   \ \\
  \end{split}\]
  (where the last equality follows from Lemma \ref{FY}).
  
  To finish the proof of the lemma, we now claim that for any (ordinary or special) Gushel--Mukai threefold $Y$ we have
    \begin{equation}\label{claim} \Delta_Y\cdot h\ \ \in \bigl\langle \ima\bigl(A^\ast(\Gr(2,5))\to A^\ast(Y)\bigr)\bigr\rangle\ .\end{equation}
    Combined with Lemma \ref{FY}, this means that for a special Gushel--Mukai threefold $Y$ (and $\YY\to B$ as above) there is equality
    \[  GDA^\ast_B(Y\times Y)= \bigl\langle h\bigr\rangle \oplus \QQ[\Delta_Y]\ .\]
    Then, since the diagonal is linearly independent in cohomology of $\bigl\langle h\bigr\rangle$ (since $h^{1,2}(Y)\not=0$), this proves the lemma.
    
    It remains to prove the claim \eqref{claim}. Using the spread argument \cite[Lemma 3.2]{Vo}, it suffices to prove equality \eqref{claim} for the very general Gushel--Mukai threefold. Thus, we may assume that $Y$ is ordinary, and moreover that
     \[ Y=Y^\prime\cap Q\ ,\]
     where $Q$ is a quadric and $Y^\prime=\Gr(2,5)\cap H_1\cap H_2$ is a smooth fourfold (where $H_1, H_2$ are Pl\"ucker hyperplanes) and $Y^\prime$ is such that
     \[A^\ast(Y^\prime)=\ima\bigl(A^\ast(\Gr(2,5))\to A^\ast(Y^\prime)\bigr)\ .\] 
     (Indeed, the smooth fourfold $Y^\prime$ has trivial Chow groups \cite[Corollary 4.6]{56}, and the very general $Y^\prime$ has no primitive cohomology, as follows from \cite[Lemma 3.15]{56}).
     The excess intersection formula then implies that
     \[ \Delta_Y\cdot h = {1\over 2}\,  \Delta_{Y^\prime}\vert_{Y\times Y}\ ,\]
     and the claim \eqref{claim} follows. 
   \end{proof}

  Lemma \ref{FY2} being proven, all conditions of Proposition \ref{crit} are met with, and so fibers $Y$ of the family $\YY\to B$ have an MCK decomposition; this settles (\rom5). 
     \end{proof}

 \section{The tautological ring}
 
 \begin{theorem}\label{taut} Let $Y$ be a Fano threefold of Picard number 1. Assume that $Y$ has an MCK decomposition, and $Y$ is member of a family $\YY\to B$  such that
 $\YY\times_B \YY\to B$ has the Franchetta property.
 For $m\in\NN$, let
  \[ R^\ast(Y^m):=\bigl\langle (p_i)^\ast(h), (p_{ij})^\ast(\Delta_Y)\bigr\rangle\ \subset\ \ \ A^\ast(Y^m)   \]
  be the $\QQ$-subalgebra generated by pullbacks of the polarization $h\in A^1(Y)$ and pullbacks of the diagonal $\Delta_Y\in A^3(Y\times Y)$. (Here $p_i$ and $p_{ij}$ denote the various projections from $Y^m$ to $Y$ resp. to $Y\times Y$).
  The cycle class map induces injections
   \[ R^\ast(Y^m)\ \hookrightarrow\ H^\ast(Y^m,\QQ)\ \ \ \hbox{for\ all\ }m\in\NN\ .\]
   \end{theorem}

\begin{proof} This is inspired by the analogous result for cubic hypersurfaces \cite[Section 2.3]{FLV3}. In its turn, the result of \cite{FLV3} was inspired by analogous results for hyperelliptic curves \cite{Ta2}, \cite{Ta} (cf. Remark \ref{tava} below) and for K3 surfaces \cite{V17}, \cite{Yin}.

Let $d$ denote the degree of $Y$, and let $2b:=\dim H^3(Y,\QQ)$. As in \cite[Section 2.3]{FLV3}, let us write $o:={1\over d} h^3\in A^3(Y)$ (the ``distinguished zero-cycle'') and
  \[ \tau:= \Delta_Y - {1\over d}\, \sum_{j=0}^3  h^j\times h^{3-j}\ \ \in\ A^3(Y\times Y) \]
  (this cycle $\tau$ is nothing but the projector on the motive $h^3(Y)$ considered above).
Moreover, let us write 
  \[ \begin{split}    h_i&:=(p_i)^\ast(h)\ \ \in \ A^1(Y^m)\ ,\\
                                o_i&:= (p_i)^\ast(o)\ \ \in\ A^3(Y^m)\ ,\\
                                               \tau_{i,j}&:=(p_{ij})^\ast(\tau)\ \ \in\ A^3(Y^m)\ .\\
                         \end{split}\]
We define the $\QQ$-subalgebra
  \[ \bar{R}^\ast(Y^m):=\langle o_i, h_i, \tau_{i,j}\rangle\ \ \ \subset\ H^\ast(Y^m,\QQ) \]
  (where $i$ ranges over $1\le i\le m$, and $1\le i<j\le m$). One can prove (just as \cite[Lemma 2.11]{FLV3} and \cite[Lemma 2.3]{Yin}) that the $\QQ$-algebra $ \bar{R}^\ast(Y^m)$
  is isomorphic to the free graded $\QQ$-algebra generated by $o_i,h_i,\tau_{ij}$, modulo the following relations:
    \begin{equation}\label{E:X'}
			o_i\cdot o_i = 0, \quad h_i \cdot o_i = 0,  \quad 
			h_i^3 =d\,o_i\,;
			\end{equation}
			\begin{equation}\label{E:X2'}
			\tau_{i,j} \cdot o_i = 0 ,\quad \tau_{i,j} \cdot h_i = 0, \quad \tau_{i,j} \cdot \tau_{i,j} = 2b\, o_i\cdot o_j
			\,;
			\end{equation}
			\begin{equation}\label{E:X3'}
			\tau_{i,j} \cdot \tau_{i,k} = \tau_{j,k} \cdot o_i\,;
			\end{equation}
			\begin{equation}\label{E:X4'}
			\sum_{\sigma \in \mathfrak{S}_{2b+2}}  \prod_{i=1}^{b+1} \tau_{\sigma(2i-1), \sigma(2i)} = 0\,. 
			\end{equation}

To prove Theorem \ref{taut}, we need to check that these relations are also verified modulo rational equivalence.
The relations \eqref{E:X'} take place in $R^\ast(Y)$ and so they follow from the Franchetta property for $Y$. 
The relations \eqref{E:X2'} take place in $R^\ast(Y^2)$. The first and the last relations are trivially verified, because $Y$ being Fano one has $A^6(Y^2)=\QQ$. As for the second relation of \eqref{E:X2'}, this follows from the Franchetta property for $Y\times Y$. (Alternatively, it is possible to deduce the second relation from the MCK decomposition: indeed, the product $\tau_{} \cdot h_i$ lies in $A^4_{(0)}(Y^2)$, and it is readily checked that $A^4_{(0)}(Y^2)$ injects into $H^8(Y^2,\QQ)$.)
   
   Relation \eqref{E:X3'} takes place in $R^\ast(Y^3)$ and follows from the MCK relation. Indeed, we have
   \[  \Delta_Y^{sm}\circ (\pi^3_Y\times\pi^3_Y)=   \pi^6_Y\circ \Delta_Y^{sm}\circ (\pi^3_Y\times\pi^3_Y)  \ \ \ \hbox{in}\ A^6(Y^3)\ ,\]
   which (using Lieberman's lemma) translates into
   \[ (\pi^3_Y\times \pi^3_Y\times\Delta_Y)_\ast    \Delta_Y^{sm}  =   ( \pi^3_Y\times \pi^3_Y\times\pi^6_Y)_\ast \Delta_Y^{sm}                                                         \ \ \ \hbox{in}\ A^6(Y^3)\ ,\]
   which means that
   \[  \tau_{1,3}\cdot \tau_{2,3}= \tau_{1,2}\cdot o_3\ \ \ \hbox{in}\ A^6(Y^3)\ .\]

 
 It is left to consider relation \eqref{E:X4'}, which takes place in $R^\ast(Y^{2b+2})$.
 To check that this relation is also verified modulo rational equivalence, 
 we observe that relation \eqref{E:X4'} involves a cycle  
 contained in
   \[ A^\ast\bigl(\sym^{2b+2} (h^3(Y)\bigr)\ .\]
 But we have vanishing of the Chow motive
   \[ \sym^{2b+2} h^3(Y)=0\ \ \hbox{in}\ \MM_{\rm rat}\ ,\] 
 because $\dim H^3(Y,\QQ)=2b$ and $h^3(Y)$ is oddly
  finite-dimensional in the sense of Kimura \cite{Kim} (all Fano threefolds are known to have Kimura finite-dimensional motive \cite[Theorem 4]{43}).  
  This establishes relation \eqref{E:X4'}, modulo rational equivalence, and ends the proof.
%
%
%
 \end{proof}

\begin{remark}\label{tava} Given a curve $C$ and an integer $m\in\NN$, one can define the {\em tautological ring\/}
  \[ R^\ast(C^m):= \bigl\langle  (p_i)^\ast(K_C),(p_{ij})^\ast(\Delta_C)\bigr\rangle\ \ \ \subset\ A^\ast(C^m) \]
  (where $p_i, p_{ij}$ denote the various projections from $C^m$ to $C$ resp. $C\times C$).
  Tavakol has proven \cite[Corollary 6.4]{Ta} that if $C$ is a hyperelliptic curve, the cycle class map induces injections
    \[  R^\ast(C^m)\ \hookrightarrow\ H^\ast(C^m,\QQ)\ \ \ \hbox{for\ all\ }m\in\NN\ .\]
   On the other hand, there are many (non hyperelliptic) curves for which the tautological ring $R^\ast(C^3)$ does {\em not\/} inject into cohomology (this is related to the non-vanishing of the Ceresa cycle, cf. \cite[Remark 4.2]{Ta} and also \cite[Example 2.3 and Remark 2.4]{FLV2}). 
\end{remark}

%
%
%
%
%

\newpage
\section{A table}

Table 1 below lists all Fano threefolds with Picard number 1 (the classification of Fano threefolds is contained in \cite{IP}). The last column indicates the existence of an MCK decomposition. Note that a Fano threefold $X$ with $h^{1,2}(X)=0$ has trivial Chow groups (i.e. $A^\ast_{hom}(X)=0$), and so these Fano threefolds have an MCK decomposition for trivial reasons. The asterisks indicate new cases settled in this paper. Question marks indicate cases I am not able to settle.

\medskip
\vskip1cm

\begin{table}[h]
\centering
\begin{tabular}{||c c c c c c||} 
 \hline
 Label & Index & Degree & $h^{1,2}$   & Description  & MCK
  \\ 
 [0.5ex] 
 \hline\hline
 4 & 4 & 1 & 0 & $\PP^3$ &trivial \\
 \hline
 3 & 3 & 2 & 0 & $X_2\subset\PP^4$ & trivial\\
 \hline
 2.1 & 2 & 1 & 21 &$X_6\subset\PP(1^3,2,3)$ & $\ast$\\
 \hline
 2.2 & 2 & 2 &10& $X_4\subset\PP(1^4,2)$ & $\ast$\\
 \hline
 2.3 & 2 & 3 &5& $X_3\subset\PP^4$ & \cite{Diaz}, \cite{FLV2}\\
 \hline
 2.4& 2 & 4 & 2& $X_{(2,2)}\subset\PP^5$ & \cite{2q}\\
 \hline
 2.5 & 2&5 & 0 &$\Gr(2,5)\cap L\subset\PP^9$& trivial\\
 \hline
 1.2 & 1&2 & 52&$X_6\subset\PP(1^4,3)$& $\ast$\\
 \hline
 1.4.a&1&4&30&$X_4\subset\PP^4$& ?\\
 \hline
 1.4.b&1&4&30&$X\xrightarrow{2:1} Q\hbox{\ with\ quartic\ branch\ locus}$&$\ast$\\
 \hline
 1.6&1&6&20&$X_{(2,3)}\subset\PP^5$ & \cite{55}\\
 \hline
 1.8&1&8&14&$X_{(2,2,2)}\subset\PP^6$ & ?\\
 \hline
 1.10.a&1&10&10&ordinary Gushel--Mukai 3fold&?\\
 \hline
  1.10.b&1&10&10&special Gushel--Mukai 3fold&$\ast$\\
  \hline
  1.12&1&12&7& $\OGr_+(5,10)\cap L\subset\PP^{15}$&?\\
  \hline
  1.14&1&14&5&$\Gr(2,6)\cap L\subset\PP^{14}$&\cite{g8}\\
  \hline
  1.16&1&16&3&$\LGr(3,6)\cap L\subset\PP^{13}$&?\\
  \hline
  1.18&1&18&2&$G_2/P\cap L\subset\PP^{13}$&\cite{g10}\\
  \hline
  1.22&1&22&0&$V(s)\subset\Gr(3,7)$&trivial\\
  [1ex] 
 \hline
\end{tabular}
\caption{All Fano threefolds with Picard number 1. Here, $X_{(d_1,\ldots,d_r)}$ denotes a complete intersection of multidegree $(d_1,\ldots,d_r)$, $Q$ is a quadric, and $L\subset\PP^r$ is a linear subspace of the appropriate dimension. The notations $\LGr(3,6)$ and $\OGr_+(5,10)$ indicate the Lagrangian Grassmannian, resp. a connected component of the orthogonal Grassmannian. In 1.22, $V(s)$ denotes the zero locus of a section of some vector bundle.}
\label{table:1}
\end{table}
\medskip    
 \vskip1cm
\begin{nonumberingt} Thanks to Mr. Kai Laterveer of the Lego University of Schiltigheim who provided inspiration for this work.
\end{nonumberingt}

\end{document}